\documentclass[11pt]{article}
\usepackage{amsfonts}
\usepackage{amsmath}
\usepackage{amssymb}
\begin{document}

\baselineskip 18pt
\def\o{\over}
\def\e{\varepsilon}
\title{\Large\bf  Shifted\ Character\ Sums\ with\ Multiplicative\
Coefficients}
\author{Ke\ Gong\ and\ Chaohua\ Jia}
\date{}
\maketitle {\small \noindent {\bf Abstract.} Let $f(n)$ be a
multiplicative function satisfying $|f(n)|\leq 1$, $q$ $(\leq N^2)$
be a prime number and $a$ be an integer with $(a,\,q)=1$, $\chi$ be
a non-principal Dirichlet character modulo $q$. In this paper, we
shall prove that
$$
\sum_{n\leq N}f(n)\chi(n+a)
\ll {N\o q^{1\o 4}}\log\log(6N)+q^{1\o 4}N^{1\o
2}\log(6N)+{N\o \sqrt{\log\log(6N)}}.
$$

We shall also prove that
\begin{align*}
&\sum_{n\leq N}f(n)\chi(n+a_1)\cdots\chi(n+a_t)\ll {N\o q^{1\o
4}}\log\log(6N)\\
&\quad+q^{1\o 4}N^{1\o 2}\log(6N)+{N\o \sqrt{\log\log(6N)}},
\end{align*}
where $t\geq 2$, $a_1,\,\ldots,\,a_t$ are pairwise distinct integers
modulo $q$. }

\vskip.3in
\noindent{\bf 1. Introduction}

Let $q$ be a prime number, $a$ be an integer with $(a,\,q)=1$,
$\chi$ be a non-principal Dirichlet character modulo $q$.

Since the 1930s, I. M. Vinogradov had begun the study on character
sums over shifted primes
$$
\sum_{p\leq N}\chi(p+a),
$$
and obtained deep results [8, 9], where $p$ runs through prime
numbers. His best known result is a nontrivial estimate for the
range $N^\e\leq q\leq N^{{4\o 3}-\e}$, where $\e$ is a sufficiently
small positive constant, which lies deeper than the direct
consequence of Generalized Riemann Hypothesis. Later, Karatsuba [4]
widen the range to $N^\e\leq q\leq N^{2-\e}$, where Burgess's method
was applied.

For the M\"obius function $\mu(n)$, one can get same results on
sums
$$
\sum_{n\leq N}\mu(n)\chi(n+a)
$$
as that on sums over shifted primes.

In this paper, we consider the general sum
$$
\sum_{n\leq N}f(n)\chi(n+a), \eqno (1.1)
$$
where $f(n)$ is a multiplicative function satisfying $|f(n)|\leq 1$.
We shall apply the method in Section 2 in [2], which is called as
the finite version of Vinogradov's inequality, to give a nontrivial
estimate for the sum (1.1) when $q$ is in a suitable range.

{\bf Theorem 1}. Let $f(n)$ be a multiplicative function satisfying
$|f(n)|\leq 1$, $q$ $(\leq N^2)$ be a prime number and $a$ be an
integer with $(a,\,q)=1$, $\chi$ be a non-principal Dirichlet
character modulo $q$. Then we have
\begin{align*}
\qquad\qquad\qquad&\sum_{n\leq N}f(n)\chi(n+a) \ll {N\o q^{1\o
4}}\log\log(6N)\qquad\qquad\qquad\qquad\ (1.2)\\
&\quad+q^{1\o 4}N^{1\o 2}\log(6N) +{N\o \sqrt{\log\log(6N)}}.
\end{align*}

{\bf Remark}. The estimate in (1.2) is nontrivial for
$$
(\log\log(6N))^{4+\e}\ll q\ll {N^2\o (\log(6N))^{4+\e}},
$$
which should be compared with the conjectural nontrivial range as
indicated by Karatsuba [4, p.\,325].

In the same way as Theorem 1, we shall prove the following

{\bf Theorem 2}. We assume that $f(n),\,q,\,\chi$ same as in Theorem
1. For $t\geq 2$ and pairwise distinct integers $a_1,\,\ldots,\,a_t$
modulo $q$, we have
\begin{align*}
\qquad\qquad&\sum_{n\leq N}f(n)\chi(n+a_1)\cdots\chi(n+a_t)\ll {N\o
q^{1\o 4}}\log\log(6N)\qquad\qquad\ \  (1.3)\\
&\quad +q^{1\o 4}N^{1\o 2}\log(6N) +{N\o \sqrt{\log\log(6N)}},
\end{align*}
where implied constant depends on $t$.

{\bf Remarks}. 1) Taking $f=\mu$ to be the M\"obius function in
Theorem 2, we obtain an example for the M\"obius Randomness Law.
Such an example answers, in a special case, a problematic issue
posed by P.~Sarnak, see [7, p.\,3].

2) The $t=2$ case of Theorem 2 corresponds to Karatsuba [5]. While
for any $t\geq 3$, our result should be compared with a conditional
result of Karatsuba [6], which relies on a conjectural upper bound
for a kind of character sums in two variables.

Throughout this paper, we assume that $N$ is sufficiently large and
set
\begin{align*}
\qquad\quad &d_0=\sqrt{\log\log(6N)},\quad
D_0=e^{d_0}=\exp(\sqrt{\log\log(6N)}),
\qquad\qquad (1.4)\\
&d_1=d_0^2=\log\log(6N),\quad D_1=e^{d_1}=\log(6N).
\end{align*}
Let $p,\,q$ denote prime numbers, $\e$ be a sufficiently small
positive constant.

\vskip.3in
\noindent{\bf 2. The proof of Theorem 1}

{\bf Lemma 1}. Let $f(n)$ be a multiplicative function satisfying
$|f(n)|\leq 1$, $q$ be a prime number, $a$ be an integer, $\chi$ be
a Dirichlet character modulo $q$. $d_0,\,d_1$ are defined as in
(1.4). Then we have
\begin{align*}
\qquad &\ \,\sum_{n\leq N}f(n)\chi(n+a)\qquad\qquad\qquad\qquad
\qquad\qquad\qquad\qquad\qquad\qquad (2.1)\\
&\ll\sum_{r=[d_0]+1}^{[d_1]-1}\sum_{y\leq {N\o e^r}}\Bigl|\sum_
{\substack{e^r\leq p< e^{r+1}\\ (p,\,q)=1\\ p\leq{N\o
y}}}f(p)\chi(py+a)\Bigr| +{N\o \sqrt{\log\log(6N)}}+{N\o q}.
\end{align*}

Proof. We have
$$
\sum_{n\leq N}f(n)\chi(n+a)=\sum_{\substack{n\leq N\\
(n,\,q)=1}}f(n)\chi(n+a)+\sum_{\substack{n\leq N\\
(n,\,q)>1}}f(n)\chi(n+a)
$$
and
$$
\sum_{\substack{n\leq
N\\ (n,\,q)>1}}f(n)\chi(n+a)\ll\sum_{\substack{n\leq N\\
q|n}}1\ll{N\o q}+1.
$$

By the discussion in Section 2 of [3], the idea of which was adopted
from [2] with some modification, we get
\begin{align*}
\sum_{\substack{n\leq N\\
(n,\,q)=1}}f(n)\chi(n+a)&\ll\sum_{r=[d_0]+1}^{[d_1]-1}\sum_{y\leq
{N\o e^r}}\Bigl|\sum_ {\substack{e^r\leq p< e^{r+1}\\ (p,\,q)=1\\
p\leq{N\o y}}}f(p)\chi(py+a)\Bigr|\\
&+{N\o \sqrt{\log\log(6N)}}.
\end{align*}
Hence, the conclusion of Lemma 1 follows.

{\bf Lemma 2}. Let $q$ be a prime number, $\chi_1,\,\ldots,\,\chi_r$
be Dirichlet characters modulo $q$, at least one of which is
non-principal. Let $f(X)\in{\Bbb F}_q[X]$ be an arbitrary polynomial
of degree $d$. Then for pairwise distinct $a_1,\,\ldots,\,a_r\in
{\Bbb F}_q$, we have
$$
\Bigl|\sum_{x\in{\Bbb F}_q}\chi_1(x+a_1)\cdots\chi_r(x+a_r)
e\Bigl({f(x)\o q}\Bigr)\Bigr|\leq (r+d)q^{1\o 2}.
$$

This is Lemma 17 in [1].

{\bf Lemma 3}. Let $q$ be a prime number, $\chi$ be a non-principal
Dirichlet character modulo $q$. Then for an arbitrary integer $h$
with $1\leq h\leq q$ and distinct $s,\,t\in{\Bbb F}_q$,
$$
\sum_{x=1}^h\chi\Bigl({x+s\o x+t}\Bigr)=O(q^{1\o 2}\log q)
$$
holds true. Here we write ${1\o n}$ as the multiplicative inverse of
$n$ such that ${1\o n}\cdot n\equiv 1\,({\rm mod}\,q)$ and appoint
${1\o 0}=0$.

This is Lemma 18 in [1].

{\bf Lemma 4}. Let $q$ be a prime number, $\chi$ be a non-principal
Dirichlet character modulo $q$, $(a,\,q)=1$. Then for two primes
$p_1,\,p_2$ with $(p_1,\,q)=(p_2,\,q)=1,\,p_1\not\equiv p_2\,({\rm
mod}\,q)$, we have
$$
\sum_{X<y\leq Z}\chi\Bigl({p_1y+a\o p_2y+a}\Bigr)\ll {Z-X\o q}
\sqrt{q}+\sqrt{q}\log q.
$$

Proof. By Lemmas 2 and 3, we have
\begin{align*}
\sum_{X<y\leq Z}\chi\Bigl({p_1y+a\o p_2y+a}\Bigr)&=\chi(p_1){\bar
\chi}(p_2)\sum_{X<y\leq Z}\chi\Bigl({y+a\overline{p_1}\o
y+a\overline{p_2}}\Bigr)\\
&\ll{Z-X\o q} \sqrt{q}+\sqrt{q}\log q.
\end{align*}

Let
$$
Y={N\o e^r}. \eqno (2.2)
$$
We shall estimate the sum
$$
{\sum}_1=\sum_{y\leq Y}\Bigl|\sum_{\substack{e^r \leq p< e^{r+1}\\
(p,\,q)=1\\ p\leq{N\o y}}}f(p)\chi(py+a)\Bigr|. \eqno (2.3)
$$

By Cauchy's inequality, we have
$$
{\sum}_1\leq Y^{1\o 2}\Bigl(\sum_{y\leq Y}\Bigl|\sum_{\substack {e^r
\leq p< e^{r+1}\\ (p,\,q)=1\\ p\leq{N\o
y}}}f(p)\chi(py+a)\Bigr|^2\Bigr)^{1\o 2}.
$$
An application of Lemma 4 to
$$
{\sum}_2=\sum_{y\leq Y}\Bigl|\sum_{\substack{e^r \leq p< e^{r+1}\\
(p,\,q)=1\\ p\leq{N\o y}}}f(p)\chi(py+a)\Bigr|^2
$$
produces
\begin{align*}
{\sum}_2&=\sum_{y\leq Y}\sum_{\substack{e^r \leq p_1< e^{r+1}\\
(p_1,\,q)=1\\ p_1\leq{N\o y}}}\sum_{\substack{e^r \leq p_2< e^{r+1}\\
(p_2,\,q)=1\\ p_2\leq{N \o y}}}f(p_1)\overline{f(p_2)}\chi\Bigl({p_1y
+a\o p_2y+a}\Bigr)\\
&=\sum_{\substack{e^r \leq p_1<e^{r+1}\\
(p_1,\,q)=1}}\sum_{\substack{e^r \leq p_2< e^{r+1}\\ (p_2,\,q)=1}}
f(p_1)\overline{f(p_2)}\sum_{\substack{y\leq Y\\
y\leq{N\o \max(p_1,\,p_2)}}}\chi\Bigl({p_1y+a\o p_2y+a}\Bigr)\\
&\ll\sum_{\substack{e^r \leq p_1< e^{r+1}\\
(p_1,\,q)=1}}\sum_{\substack{e^r \leq p_2< e^{r+1}\\
(p_2,\,q)=1}}\Bigl| \sum_{\substack{y\leq Y\\ y\leq{N\o
\max(p_1,\,p_2)}}}\chi\Bigl(
{p_1y+a\o p_2y+a}\Bigr)\Bigr|\\
&\ll\sum_{e^r \leq p_1< e^{r+1}}\sum_{\substack{e^r \leq p_2<
e^{r+1}\\ p_2\equiv p_1\,({\rm mod}\,q)}}Y\\
&\ +\sum_{\substack{e^r \leq p_1< e^{r+1}\\
(p_1,\,q)=1}}\sum_{\substack{e^r \leq p_2< e^{r+1}\\
(p_2,\,q)=1\\ p_2\not\equiv p_1\,({\rm
mod}\,q)}}\Bigl|\sum_{\substack {y\leq Y\\ y\leq{N\o
\max(p_1,\,p_2)}}}\chi\Bigl({p_1y+a\o p_2y+
a}\Bigr)\Bigr|\\
&\ll Ye^r\Bigl({e^r\o q}+1\Bigr)+\sum_{e^r \leq p_1<e^{r+1}}
\sum_{e^r\leq p_2< e^{r+1}}\Bigl({Y\o \sqrt{q}}+\sqrt{q}\log q
\Bigr)\\
&\ll Ye^r+{Ye^{2r}\o \sqrt{q}}+\sqrt{q}(\log q)e^{2r}.
\end{align*}
Thus
\begin{align*}
{\sum}_1&\ll Y^{1\o 2}\Bigl({Ye^{2r}\o \sqrt{q}}+Ye^r+\sqrt{q}
(\log q)e^{2r}\Bigr)^{1\o 2}\\
&\ll {Ye^r\o q^{1\o 4}}+Ye^{r\o 2}+q^{1\o 4}(\log q)^{1\o 2}Y^{1\o
2}e^r\\
&\ll{N\o q^{1\o 4}}+{N\o e^{r\o 2}}+q^{1\o 4}(\log q)^{1\o 2}N^{1\o
2}e^{r\o 2}.
\end{align*}
Applying this estimate to (2.1), we get
\begin{align*}
&\ \,\sum_{n\leq N}f(n)\chi(n+a)\\
&\ll\sum_{r=[d_0]+1}^{[d_1]-1}\Bigl({N\o q^{1\o 4}}+{N\o e^{r\o 2}}
+q^{1\o 4}(\log q)^{1\o 2}N^{1\o 2}e^{r\o 2}\Bigr)+{N\o
\sqrt{\log\log(6N)}}+{N\o q}\\
&\ll{N\o q^{1\o 4}}\log\log(6N)+{N\o \exp{({1\o
2}\sqrt{\log\log(6N)})}}\\
&\ +q^{1\o 4}(\log q)^{1\o 2}N^{1\o 2}(\log(6N))^{1\o 2}
+{N\o \sqrt{\log\log(6N)}}\\
&\ll{N\o q^{1\o 4}}\log\log(6N)+q^{1\o 4}N^{1\o 2}\log(6N)
+{N\o \sqrt{\log\log(6N)}}.
\end{align*}

So far the proof of Theorem 1 is complete.

\vskip.3in
\noindent{\bf 3. The proof of Theorem 2}

{\bf Lemma 5}. Let $q$ be a prime number, $\chi$ be a non-principal
Dirichlet character modulo $q$. Then for an arbitrary integer $h$
with $1\leq h\leq q$ and pairwise distinct
$a_1,\,\ldots,\,a_t,\,b_1,\,\ldots,\,b_t(t\geq 2)\in{\Bbb F}_q$,
$$
\sum_{x=1}^h\chi\Bigl({(x+a_1)\cdots(x+a_t)\o
(x+b_1)\cdots(x+b_t)}\Bigr)=O(q^{1\o 2}\log q)
$$
holds true, where the implied constant depends on $t$.

Proof. It comes from Lemma 2 in the same way as Lemma 3.

{\bf Lemma 6}. Let $q$ be a prime number, $\chi$ be a non-principal
Dirichlet character modulo $q$. Assume that $t\geq
2,\,(a_1,\,q)=\cdots=(a_t,\,q)=1$, $a_1,\,\ldots,\,a_t$ are pairwise
distinct modulo $q$ and that primes $p_1,\,p_2$ satisfy
$(p_1,\,q)=(p_2,\,q)=1$, $p_1\not\equiv p_2\,({\rm mod}\,q)$. Then if
$p_2\not\equiv \overline{a_i}a_j p_1\,({\rm mod}\,q)$ $(1\leq i,\,j\leq
t)$, we have
$$
\sum_{X<y\leq Z}\chi\Bigl({(p_1y+a_1)\cdots(p_1y+a_t)\o (p_2y
+a_1)\cdots(p_2y+a_t)}\Bigr)\ll {Z-X\o q} \sqrt{q}+\sqrt{q}\log q.
$$
If $p_2\equiv \overline{a_i}a_j p_1\,({\rm mod}\,q)$ for some
$(i,\,j)$, we have trivial bound
$$
\sum_{X<y\leq Z}\chi\Bigl({(p_1y+a_1)\cdots(p_1y+a_t)\o (p_2y
+a_1)\cdots(p_2y+a_t)}\Bigr)\ll Z-X+1.
$$

Proof. It can be proved in the same way as Lemma 4.

We can prove Theorem 2 in the same way as Theorem 1, by discussing
the case in which one of $a_i$ is equivalent to $0\,({\rm mod}\,q)$
and the case in which no $a_i$ is equivalent to $0\,({\rm mod}\,q)$
respectively.

\vskip.3in
\noindent{\bf Acknowledgements}

The first author is supported by the National Natural Science
Foundation of China (Grant No. 11126150) and the Natural Science
Foundation of the Education Department of Henan Province (Grant No.
2011A110003). The second author is supported by the National Key
Basic Research Program of China (Project No. 2013CB834202) and the
National Natural Science Foundation of China (Grant No. 11371344 and
Grant No. 11321101).


\bigskip

\

Ke Gong

Department of Mathematics, Henan University, Kaifeng, Henan 475004,
P. R. China

E-mail: kg@henu.edu.cn

\

Chaohua Jia

Institute of Mathematics, Academia Sinica, Beijing 100190, P. R.
China

Hua Loo-Keng Key Laboratory of Mathematics, Chinese Academy of
Sciences, Beijing 100190, P. R. China

E-mail: jiach@math.ac.cn

\end{document}